\documentclass[12pt,leqno]{article}
\usepackage[official]{eurosym}
\usepackage{fullpage,amsmath,theorem}
\usepackage{epic,eepic}
\usepackage{amssymb,amsmath}
\usepackage{cases}
\def\llvdash{{\|\hskip-2pt \raise 3pt\hbox{\vrule
height 0.25pt width 0.4cm}}}

\def\bet{{\beta}}

\def\l{\langle}
\def\r{\rangle}

\def\l{{\langle}}
\def\r{{\rangle}}

\def\oa{{\overline A^{\,\lower 7pt_{\hbox{$\scriptstyle\bet}}
\hbox{$\scriptstyle 0\tau$}}}}

\def\bet{\beta}

\def\llvdash{{\|\hskip-2pt \raise 3pt\hbox{\vrule height 0.25pt
width 0.4cm}}}

\newtheorem{theorem}{Theorem}[section]

{\theorembodyfont{\rmfamily}
\newtheorem{fact}[theorem]{Fact}
{\theorembodyfont{\rmfamily}
\newtheorem{definition}[theorem]{Definition}}
{\theorembodyfont{\rmfamily}
}
{\theorembodyfont{\rmfamily}
\newtheorem{claim}{Claim}}
{\theorembodyfont{\rmfamily}
}
{\theorembodyfont{\rmfamily}
}
 \DeclareMathOperator{\dom}{dom}
\DeclareMathOperator{\rng}{rng}

\newcommand{\pr}{\medskip\noindent\textit{Proof}. }

\newcommand{\lusim}[1]{\smash{\underset{\raisebox{1.2pt}[0cm][0cm]{$\sim$}}
{{#1}}}}

\def\dom{{\rm dom}}

\def\rng{{\rm rng}}
\def\Ult{{\rm Ult}}

\def\llvdash{{\|\hskip-2pt \raise 3pt\hbox{\vrule
height 0.25pt width 0.15cm}}}

\def\Vdashbks{\hbox{$\Vdash\!\!\!\!{\raise2pt\hbox
{$\scriptscriptstyle\backslash$}}$}}

\overfullrule=0pt
\begin{document}

\title{On ultrafilters in ZF models and indecomposable ultrafilters }

\baselineskip=18pt
\author{ Eilon Bilinsky and Moti Gitik\footnote{ The work was partially supported  by ISF grant No. 882/22.
} }

\date{\today}
\maketitle

\begin{abstract}
We use indecomposable ultrafilters
to answer some questions from Hayut, Karagila \cite{H-K}.  It is shown that the bound on the strength of Usuba \cite{U} is optimal.
\end{abstract}

\section{ On indecomposable ultrafilters }

In sixties C. Chang and J. Keisler formulated the following notions:

\begin{definition}
Let $U$ be an ultrafilter on a set $I$.
\begin{enumerate}
  \item $U$ is called $(\kappa, \lambda)-$regular iff there is subset of $U$ of cardinality $\lambda$ such that any $\kappa-$members of it have empty intersection.
  \item $U$ is called $\lambda-$descendingly incomplete iff there are $\{X_\alpha \mid \alpha<\lambda\}\subseteq U$
  such that $\alpha<\beta \rightarrow X_\alpha\supseteq X_\beta$ and $\bigcap_{\alpha<\lambda}X_\alpha=\emptyset$.
  \item $U$ is $\lambda-$decomposable iff there is a partition of $I$ into disjoint sets $\l I_\alpha \mid \alpha<\lambda \r$, so that whenever $S\subseteq \lambda$ and $|S|<\lambda$, $\bigcup_{\alpha\in S}I_\alpha \not \in U$.

      \item Suppose $\delta<\lambda$ are cardinals. $ U$  is called $(\delta, \lambda)-$indecomposable if any partition \\$\l I_\nu \mid \nu<\alpha\r$ of $I$ with $\alpha<\lambda$ has a subsequence $\l I_{\nu_\xi}\mid  \xi<\beta\r$ with
$\beta<\delta$ whose union belongs to $U$.
\end{enumerate}

\end{definition}

Let state some known  facts which are relevant for us here:

\begin{fact}
 $U$ is $\lambda-$decomposable, then $U$ is $\lambda-$descendingly incomplete.
 \\If $\lambda$ is regular, then the converse holds as well.

\end{fact}

\begin{fact}\label{prop-r-k}
An ultrafilter $U$ over $I$ is $\lambda-$decomposable iff it Rudin-Keisler above a uniform ultrafilter over $\lambda$.

\end{fact}

\begin{fact}\label{prop-reg-dec}
If $U$ is $(\kappa,\lambda)-$regular ultrafilter and $\nu$ is a regular cardinal so that $\kappa\leq \nu\leq \lambda$, then $U$ is
$\nu-$descendingly incomplete, and so, $\nu-$decompossible.

\end{fact}

\begin{fact}
\begin{enumerate}
  \item $U$ is $\gamma-$indecomposable if and only if $U$ is $(\gamma, \gamma^+)-$indecomposable.
  \item   $U$ is $(\delta, \lambda)-$indecomposable if and only if $U$ is not $\gamma-$decomposable for any cardinal $\gamma$ such that $\delta\leq \gamma<\lambda$.
\end{enumerate}

\end{fact}

We will relay on the following theorem of J. Silver:

\begin{theorem}\label{thm-silver}
Let $\delta$ and $\kappa$ be cardinals with $2^\delta<\kappa$. Suppose that $U$ is\\ a $(\delta,\kappa)-$indecomposable ultrafilter over a set $I$. Then $j_U = j_W^{M_D}\circ
j_D$ where $D$ is an ultrafilter over
a cardinal less than $\delta$ and $W$ is an $M_D-\kappa-$complete $M_D-$ultrafilter over $j_U(I)$.

\end{theorem}

Let us show the following:

\begin{theorem}\label{thm1}
Let $\delta$ and $\kappa$ be cardinals with $2^\delta<\kappa$. Suppose that $U$ is\\ a $(\delta,\kappa)-$indecomposable ultrafilter over a set $I$.
Let $P$ be a $\delta-$closed forcing such that there is $\rho<\kappa, \rho^{<\delta}=\rho,|P|\leq \rho$.
Let $G\subseteq P$ be a generic.\\ Then, in $V[G]$,  $U^*=\{A\subseteq I\mid \exists B\in U (B\subseteq A)\}$ is a $(\delta,\kappa)-$indecomposable ultrafilter over a set $I$.

\end{theorem}
\pr
Consider $j_U:V\to M_U=\Ult(V, U)$. Note that $M_U$ may be ill-founded.
By the Silver theorem, $j_U=j_W^{M_D}\circ
j_D$ where $D$ is an ultrafilter over
a cardinal less $\eta$ than $\delta$ and $W$ is an $M_D-\kappa-$complete $M_D-$ultrafilter over $j_U(I)$.
\begin{claim} $j_U{}''G$ generates a generic subset of $j_U(P)$ over $M_U$.
\end{claim}
\pr
Let $E\subseteq j_U(P)$ be a dense open subset in $M_U$. \\
We have $|P|\leq \rho$, so without loss of generality, assume that $P\subseteq \rho$. The assumption $\rho^{<\delta}=\rho$ implies that $j_D(\rho)<\kappa$.
In particular, $j_U(P)=j_D(P)$.
\\
Pick $f_E:I \to V$ which represents $E$. We can assume that $f_E$ depends on $\eta$ only.
Then we will have $|\rng(f_E)|\leq \eta.$ Use $\delta-$closure of $P$ and find $E^*\subseteq P$ which is a dense open and is contained in each dense open subset of $P$ in $\rng(f_E)$. Then, $j_U(E^*)\subseteq E$.
\\
$\square$ of the claim.

Now, exactly as in a well-founded case the elementary embeddings extend. Denote extensions by $j^*, j_D^*, j_W^*$.
\\We have, for every $X\subseteq I$,
$$X\in U \Leftrightarrow [id]_U \in j_U(X).$$
Now, in $V[G]$, let $A\subseteq I$ and $\lusim{A}$ be its name.
\\Set
$$A\in U^* \Leftrightarrow \exists p\in G (j_U(p) \Vdash[id]_U \in j_U(\lusim{A})).$$
Then, $U^*\supseteq U$.

\begin{claim}
 $U^*=\{A\subseteq I\mid \exists B\in U (B\subseteq A)\}$.
\end{claim}
\pr
Let $A\in U^*$. Pick $p \in G$ such that $j_U(p) \Vdash[id]_U \in j_U(\lusim{A})$.
\\Set $B=\{\nu\in I \mid p\Vdash \nu\in \lusim{A}\}$.
Then $B\in U$ and $B\subseteq A$, since $p\in G$.
\\
$\square$ of the claim.

The next claim completes the proof.

\begin{claim}
 $U^*$ is a $(\delta,\kappa)-$indecomposable ultrafilter in $V[G]$.
\end{claim}
\pr
Let $\l I_\nu \mid \nu<\alpha\r$ be a partition of $I$ with $\alpha<\kappa$.
We need to show that there is  a subsequence $\l I_{\nu_\xi}\mid  \xi<\beta\r$ with
$\beta<\delta$ whose union belongs to $U^*$.
\\Apply $j^*_D$ to the partition. Let $\l I'_\nu \mid \nu<j_D(\alpha)\r$ be the result. Note that $\alpha':=j_D(\alpha)<\kappa$, and so, the further embedding $j_{W^*}$ will not move $\alpha'$. Let $j_{W^*}(\l I'_\nu \mid \nu<\alpha'\r)=\l I''_\nu \mid \nu<\alpha'\r$. Then, for every $\nu<\alpha'$, $ I''_\nu =j_{W^*}(A'_\nu)$.
There must be some $\nu^*<\alpha'$ such that $[id]_U\in I''_{\nu^*}$.
\\Let $[f]_{D}$ be a function that represents $\nu^*$. We can assume that $\rng(f) \subseteq \{ I_\nu \mid \nu<\alpha\}$.
\\Then $\bigcup_{\nu\in \rng(f)}I_\nu \in U^*$. We are done, since $|\rng(f)|\leq |\dom(f)|<\delta$.
\\
$\square$ of the claim.
\\$\square$

A similar, and a simpler argument  gives the following:

\begin{theorem}\label{thm2}

Let $\delta$ and $\kappa$ be cardinals with $2^\delta<\kappa$. Suppose that $U$ is\\ a $(\delta,\kappa)-$indecomposable ultrafilter over a set $I$.
Let $P$ be a  forcing of cardinality less than the critical point of $j_U$.
Let $G\subseteq P$ be a generic.\\ Then, in $V[G]$,  $U^*=\{A\subseteq I\mid \exists B\in U (B\subseteq A)\}$ is a $(\delta,\kappa)-$indecomposable ultrafilter over a set $I$.

\end{theorem}

\section{ On existence of indecomposable ultrafilters on non-measurable cardinals }

Clearly, if $U$ is a $\kappa-$complete uniform ultrafilter over $\kappa$, then $U$ is $\lambda-$indecomposable for every $\lambda<\kappa$.
By D. Donder \cite{D}, if $\delta<\lambda$ and $\lambda$ carries a $\delta-$indecomposable uniform ultrafilter, then there is an inner model of a measurable cardinal.
\\However, a cardinal which carries such ultrafilters need not be a measurable or even large.
Thus,
S. Ben David and M. Magidor \cite{B-M} used a supercompact to construct a model in which there is a uniform ultrafilter over $\aleph_{\omega+1}$ which is $(\omega_1, \aleph_\omega)-$indecomposable.
\\H. Woodin, starting with a measurable and building on similar ideas constructed a GCH model in which there is a uniform ultrafilter over $\aleph_\omega$ which is $(\omega_1, \aleph_\omega)-$indecomposable.
\\In \cite{G}, starting with two measurables a GCH model in which there are  regular cardinals $\omega<\kappa<\lambda$ such that $\lambda$ is not measurable and carries a uniform $\sigma-$complete ultrafilter which is $\kappa-$indecomposable.

Let us briefly recall Woodin's construction and those of \cite{G}.

\textbf{A sketch of Woodin's construction}.

Start with a measurable cardinal $\kappa$. Let $F$ be a normal ultrafilter over $\kappa$. Force with the Prikry forcing with $F$.\footnote{ Actually, Woodin combines this with collapses in order to turn $\kappa$ into $\aleph_\omega$, which is unneeded for our purposes.}
Let $\l \kappa_n \mid n<\omega\r$ be a resulting Prikry sequence.
Let $V_1=V[\l \kappa_{2n} \mid n<\omega\r]$. Define a $V_1-$filter $U_1$ over $\kappa$ by setting
$$X\in U_1 \Leftrightarrow \exists n_0<\omega (X\supseteq \{\kappa_{2n+1} \mid n_0\leq n<\omega\}).$$
Using the homogeneity of the Prikry forcing, it is possible to argue that $U_1\in V_1$.
\\Pick a non-principal ultrafilter $D$ on $\omega$ and define $U$:
$$X\in U \Leftrightarrow \exists A\in D (X\supseteq \{\kappa_{2n+1} \mid n\in A\}).$$
Again, such $U$ will be in $V_1$ and it will be a uniform $(\omega_1, \kappa)-$indecomposable there.

\textbf{A sketch of the constructioin from \cite{G}}.

Start with two measurable cardinals $\kappa<\lambda$.
\\Fix normal ultrafilters $U_\kappa,U_\lambda$ over $\kappa$ and $\lambda$ respectively.
The final ultrafilter $U$ will extend $U_\kappa \times U_\lambda$.
In order to destroy measurability of $\lambda$ a type of forcing adding Suslin trees  is iterated below $\lambda$ and at $\lambda$ itself.
Below $\lambda$ branches are added to such Suslin trees and nothing is done over $\lambda$(in K. Kunen fashion).
The iteration is arranged in a special way which allows to extend the embedding $j_{ U_\kappa \times U_\lambda}$.
This will give a uniform $(\kappa^+, \lambda)-$indecomposable ultrafilter over a non-measurable cardinal $\lambda$.
In addition, such ultrafilter will be $\kappa-$complete.

Let us conclude this section with a simple construction of a uniform \\$(\omega_1, \kappa)-$indecomposable ultrafilter over a cardinal $\kappa$ of countable cofinality, and so, not measurable.

\textbf{A simple construction.}

Assume that $\kappa$ is a limit of an increasing sequence of measurable cardinals $\l \l \kappa_n \mid n<\omega\r$.
Fix a non-principal ultrafilter $D$ over $\omega$ and let $U_n$ be a normal ultrafilter over $\kappa_n$, for every $n<\omega$.
Define $U$ over $\kappa$ as follows:
$$X\in U \Leftrightarrow \{n<\omega\mid X\cap\kappa_{n}\in U_n\}\in D.$$
It is not hard to see that such $U$ will be  a uniform $(\omega_1, \kappa)-$indecomposable ultrafilter.

\section{ Applications to ZF models}

Y. Hayut and A. Karagila, in \cite{H-K}, introduced and studied the class $\mathcal{U}$ of all infinite cardinals which carry a uniform ultrafilter in ZF context.
They asked whether the following:
\\\emph{Is it possible that to have a situation when some cardinal $\kappa$ does not carry a uniform ultrafilter, $\kappa^+$ does, but $\kappa^+$ is not measurable, and is this possible without using large cardinals? In particular, is it possible that $\aleph_0$ is the only measurable cardinal, while $\aleph_1\not \in \mathcal{U}$ and
$\aleph_2 \in \mathcal{U}$?}

T. Usuba \cite{U} showed that large cardinals are needed. Namely, he proved the following:

\begin{theorem}\label{thm-U}(ZF)
If there are cardinals $\kappa<\lambda$ with $\kappa\not \in \mathcal{U}$ and $\lambda \in \mathcal{U}$, then there is an inner model with a measurable cardinal.

\end{theorem}

Usuba argued that in an inner ZFC model there is a uniform ultrafilter over $\lambda$ which is $\kappa-$indecomposible.
Then by Donder \cite{D}, there exists an inner model with a measurable cardinal.

Our aim here will be to use indecomposable ultrafilters from ZFC models in order to provide affirmative answers to remaining parts of the above question, and also, to argue that
it is impossible to improve Usuba's lower bound.

\textbf{A model in which $\aleph_0$ is the only measurable cardinal, while $\aleph_1\not \in \mathcal{U}$ and
$\aleph_2 \in \mathcal{U}$}.

Let $U$ be a uniform $(\omega_1, \kappa)-$indecomposable ultrafilter over a singular cardinal $\kappa$ of cofinality $\omega$, in a ZFC model.
\\Use symmetric extensions with collapses in a standard fashion,  in order to turn $\kappa$ into $\aleph_2$ by collapsing a cofinal in $\kappa$ sequence to $\omega_1$.
By symmetry and \ref{thm1}, $U$ will generate an ultrafilter in the extension. In addition, by standard arguments $\omega_1$ will not carry an ultrafilter.
\\If $V$ does not have measurable cardinals, then same will hold in such symmetric extension.

\textbf{A model in which $\aleph_1$ is the  measurable cardinal, while $\aleph_2\not \in \mathcal{U}$ and
$\aleph_3 \in \mathcal{U}$, \\also it carries a $\sigma-$complete ultrafilter}.

Start with two measurable cardinals $\kappa<\lambda$.
Use \cite{G} to construct a model with
a uniform $(\kappa^+, \lambda)-$indecomposable $\kappa-$complete ultrafilter $U$ over a non-measurable cardinal $\lambda$.
\\Use symmetric extensions with collapses in a standard fashion,  in order to turn $\kappa$ into $\aleph_1$ and $\lambda$ into  $\omega_3$.
By symmetry and \ref{thm1},\ref{thm2}, $U$ will generate a $\sigma-$complete ultrafilter in the extension. In addition, by standard arguments $\omega_2$ will not carry an ultrafilter.
\\Note that $\lambda=\aleph_3$ will remain non-measurable, since the collapses (their supports) used has small cardinality, and so, cannot add a branch to a Suslin tree over $\lambda$.

\end{document}